\documentclass[11pt]{amsart}
\usepackage{amssymb}
\usepackage{graphicx}
\usepackage{amsmath}
\usepackage{epsfig}
\usepackage[all]{xy}

\newtheorem{theorem}{Theorem}[section]

\newtheorem{proposition}[theorem]{Proposition}

\begin{document}

\title{}

\Large

 \centerline{\textbf{ Conjugacy problem for subgroups with}}
\centerline{\textbf{ applications to Artin groups and}}
\centerline{\textbf{braid type groups}}

\normalsize

\author{Nuno Franco}

\begin{abstract}

Let $G$ be a group endowed with a solution to the conjugacy
problem and with an algorithm which computes the centralizer in
$G$ of any element of $G$. Let $H$ be a subgroup of $G$. We give
some conditions on $H$, under which we provide a solution to the
conjugacy problem in $H$. We apply our results to some Artin
groups and braid type groups. In particular, we give explicit
solutions to the conjugacy problem in the Artin groups of type
$\tilde{A}_n$ and $\tilde{C}_n$.

\end{abstract}

\maketitle

\centerline{\today}

\bigskip


\pagestyle{headings}

\sloppy

\section{Introduction}

Let $G$ be a group. A \textit{solution for the conjugacy problem}
in $G$ consists on an algorithm which, for given $a,b\in G,$
determines whether there exists $c\in G$ such that $a=c^{-1}bc$.
An \textit{explicit solution for the conjugacy problem} is a
solution for the conjugacy problem that also determines an element
$c\in G$ such that $a=c^{-1}bc.$ Define the  \textit{centralizer},
$C_G(a)=C(a)$, of an element $a\in G$ as the set of elements in
$G$ that commute with $a$. A solution for the \textit{membership
problem} in $G$ is an algorithm that, for a given finitely
generated subgroup $H<G$ and a given $g\in G$, decides whether
$g\in H.$

Throughout this paper $G$ will denote a group endowed with an
explicit solution for the conjugacy problem and with an algorithm
that, for any element $a\in G$, computes a  generating set for the
centralizer of $a$ in $G$ (which we always assume to be finitely
generated).

The present paper concerns the following question: under which
conditions,  for a given $H<G$, we can find an explicit solution
to the conjugacy problem in $H$.

Let $H<G$ be a subgroup. We say that $H$ verifies Condition
\textbf{PC} if there exists a group $K$, a subgroup $K'$ of $K$
and an homomorphism $\phi:G\longrightarrow K$, such that:

\begin{itemize}
    \item $H=\phi^{-1}(K'),$
    \item there is a solution for the membership problem in $K$, and
    \item the subgroup $K'$ is finite.
\end{itemize}






Let $H$ be a subgroup of $G$ that verifies Condition \textbf{PC}.
Our goal is to compute an explicit solution to the conjugacy
problem in $H$, using the triple $(K,K',\phi)$ of the definition.

Remark that Condition \textbf{PC} is essentially a condition on
the group $K$ and not on the subgroup $H$. This condition on $K$
is very restrictive since there are few known classes of groups in
which we can give  a solution to the membership problem. It is
known that the following classes of groups have a solution to the
membership problem:

\begin{itemize}
    \item finite groups,
    \item free groups,
    \item abelian groups,
    \item recursively presented metabelian groups,
    \item finitely generated fully residually free groups.
\end{itemize}

The existence of a solution to the conjugacy problem is also
closely  related to the group based cryptographic systems. In
\cite{A-A-G} and \cite{K-L} are introduced some cryptographic
systems, called group based cryptographic system, in which the
encryption security is based on:

\begin{itemize}
    \item the complexity of a solution to the conjugacy problem, and
    \item the existence of a fast solution to the word problem.
\end{itemize}

A subgroup $H$ for which it does not exists a triple that verifies
Condition \textbf{PC} will assure us that it is not possible to
use our approach to attack an $H$-based cryptographic system. Even
if we use a subgroup $H$ of $G$ that verifies Condition
\textbf{PC} to implement a cryptographic system, we still have, in
general, some security advantages when comparing to an $G$-based
cryptographic system:

\begin{itemize}
    \item the same fast solution for the word problem in $H$ as in
    the group $G$, and
    \item a higher complexity of the solution of the conjugacy problem
    in $H$.
\end{itemize}

In Section 2 we present some examples of groups and subgroups in
which we can apply or not our theory. Sections 3  contains the
main results and algorithms regarding the conjugacy problem.
Finally, in Section 4, we make some brief remarks about the
complexity of the algorithms presented in this paper.

\section{Applications}

In this section we will give some examples of groups in which we
can apply our theory. These examples will be subgroups of Artin
groups, and other groups that admit braid pictures. We start with
some definitions.

Define a \textit{Coxeter matrix} of rank $n$ as a $n\times n$
symmetric matrix, $M=(m_{i,j})$, that verifies: $m_{i,i}=1$ for
all $i=1,\ldots,n,$ and $m_{i,j}\in \{2,3,\ldots\}\cup\{\infty\}$,
for all $i,j\in\{1,\ldots,n\}$, $i\neq j$. Let $M$ be a Coxeter
matrix. The \textit{Coxeter graph} associated to $M$ is the
labelled graph, $\Gamma$, defined as follows. The set of vertices
of $\Gamma$ is $\{1,\ldots,n\}$. If $m_{i,j}=2$ then there is no
edge between $i$ and $j$, if $m_{i,j}=3,$ then there is a
non-labelled edge between $i$ and $j$, and, finally, if
$m_{i,j}>3$ or $m_{i,j}=\infty$, then there is an edge between $i$
and $j$ labelled by $m_{i,j}$.

Let $G$ be a group. For $a,b\in G$ and $n\in \mathbb{N}$ we define
the word
$$prod_n(a,b)\left\{
\begin{array}{ll}
(ab)^{\frac{n}{2}} &\text{ if }n\equiv0(\text{mod }2) \\
(ab)^{\frac{n-1}{2}}a &\text{ if } n\equiv1(\text{mod }2)
\end{array}
\right.$$ Define the \textit{Coxeter group} $W$ associated to the
Coxeter graph  $\Gamma$ as the group presented by
$$W=W(\Gamma)=\left\langle
\begin{array}{cc}
s _{1},s_{2},\ldots ,s_{n} & \left|
\begin{array}{l}
s_i^2=1,\quad i\in \{1,\ldots,n\} \\
prod_{m_{i,j}}(s_{i},s_{j})=prod_{m_{j,i}}(s_{j},s_{i}) \\
\quad \quad \quad \quad \quad \quad \quad i\neq j\text{ and
}m_{i,j}\neq\infty
\end{array}
\right.
\end{array}
\right\rangle,$$ where $M=(m_{i,j})$ is the Coxeter matrix of
$\Gamma$. Define the \textit{Artin group} associated to $\Gamma$
to be the group presented by
$$A=A(\Gamma)=\left\langle
\begin{array}{cc}
s _{1},s_{2},\ldots ,s_{n} & \left|
\begin{array}{l}
prod_{m_{i,j}}(s_{i},s_{j})=prod_{m_{j,i}}(s_{j},s_{i}) \\
\quad \quad \quad \quad \quad \quad \quad \quad i\neq j\text{ and
}m_{i,j}\neq\infty
\end{array}
\right.
\end{array}
\right\rangle.$$

If the group $W(\Gamma)$ is finite, then we say that $A(\Gamma)$
is an \textit{Artin group of spherical type}. It is known that the
Artin groups of spherical type have explicit solutions for the
conjugacy problem (see \cite{B-K-L}, \cite{FGM}, \cite{G1},
\cite{Ge} and \cite{Pi2}). Moreover, we also know, for each of
these groups, algorithms that compute a finite generating set for
the centralizer of any element of $A(\Gamma)$ (see \cite{FGM2}).
So Artin groups of spherical type are examples of groups in which
we can apply our theory. A complete classification of these groups
is given in \cite{B}. Troughout this paper we will denote by
$\mu:A(\Gamma)\longrightarrow W(\Gamma)$  the canonical
epimorphism.

\begin{figure}[ht]
  \includegraphics{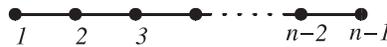}\\
  \caption{Coxeter graph $A_{n-1}$}\label{Fig1}
\end{figure}

An important example of a finite Coxeter group is the $n$th
symmetric group, $\Sigma_n$. The group $\Sigma_n$ is generated by
the simple transpositions $\tau_i=(i,i+1)$, $i=1,\ldots,n-1$. The
Coxeter relations in $\Sigma_n$ are defined by the $(n-1)\times
(n-1)$ Coxeter matrix $M=(m_{i,j})$, where $m_{i,i}=1$ for all
$i=1,\ldots,n-1$, $m_{i+1,i}=m_{i,i+1}=3$ for all
$i=1,\ldots,n-2$, and $m_{i,j}=2$ for all $i,j=1,\ldots,n-1$,
$|i-j|\geq 2$. The corresponding Coxeter graph is the graph
$A_{n-1}$ pictured in Figure \ref{Fig1} (see \cite{B}).

\vspace{0.5cm}

The Artin group associated to $\Sigma_n$ is the braid group,
$\mathcal{B}_n$, on $n$ strands defined by Artin in \cite{Artin}.
Artin's presentation for the braid group $\mathcal{B}_n$ is
$$\mathcal{B}_n=\left\langle
\begin{array}{cc}
\sigma _{1},\sigma _{2},\ldots ,\sigma _{n-1} & \left|
\begin{array}{ll}
\sigma _{i}\sigma _{j}=\sigma _{j}\sigma _{i} & (|i-j|\geq 2) \\
\sigma _{i}\sigma _{i+1}\sigma _{i}=\sigma _{i+1}\sigma _{i}\sigma
_{i+1} & (i=1,\ldots ,n-2)
\end{array}
\right.
\end{array}
\right\rangle. \text{ }  \label{presen}
$$

The braid group $\mathcal{B}_n$ admits also a geometric
interpretation as follows. Consider a disc  $D$, and  $n$ distinct
points $P_1, \ldots, P_n$ in the interior of $D$. A \textit{braid}
$b$ is the union of $n$ disjoint paths $b_1,\ldots, b_n$ in the
cylinder $D\times [0,1]$ such that, for all $i=1,\ldots,n$, the
path $b_i$ (the $i$th strand) goes monotonically in $t\in [0,1]$
from $(P_i,0)$ to some $(P_j,1)$. We consider braids modulo
isotopy, that is, two braids are equivalent if we can transform
one into the other by deforming the strands, but keeping the ends
fixed. Then the braid group is  the group  of  isotopy classes of
braids.

Note that the canonical epimorphism
$\mu:\mathcal{B}_n\longrightarrow \Sigma_n$ is defined by:
$\mu(\sigma_i)=\tau_i.$ Let $b\in \mathcal{B}_n$. We say that $b$
is a \textit{pure braid} if $\mu(b)=1$, or, equivalently, if the
$i$th strand goes form $(P_i,0)$ to $(P_i,1)$, for all
$i=1,\ldots,n$.

\bigskip

\noindent \textbf{Example 1}: Let $n\in\mathbb{N}$,
$X\subset\{1,\ldots,n\}$, and $\Sigma_{n}$ be the $n$th symmetric
group. Define $\Sigma_n(X)$ to be the subgroup of $\omega \in
\Sigma_n$ such that, for all $i\in X$, $\omega(i)=i$. Define
$\mathcal{B}_n(X)=\mu^{-1}(\Sigma_n(X))$. Briefly speaking, the
group $\mathcal{B}_n(X)$ is the set of elements in $\mathcal{B}_n$
that are pure in the $x$th strand, for all $x\in X$.

Since the group $\Sigma_n$ is finite, the triple
$(\Sigma_{n},\Sigma_n(X),\mu)$ verifies Condition \textbf{PC},
hence we can apply our theory to the subgroup $\mathcal{B}_n(X)<
\mathcal{B}_n$, and shall obtain by this way a solution to the
explicit conjugacy problem in $\mathcal{B}_n(X)$. Remark that
$\mathcal{B}_n(X)$ is a finite index subgroup of $\mathcal{B}_n$,
hence we can also give an algorithm to compute a finite generating
set of the centralizer of any element of $\mathcal{B}_n(X)$. This
example will be used later.

\bigskip

\noindent \textbf{Example 2}: The above construction has an easy
extension to all spherical type Artin groups as follows. Pick a
subgroup $K'$ of $W$. Then, as $W$ is finite, the triple
$(W,K',\mu)$ verifies Condition \textbf{PC}, and we can again
apply our theory to the subgroup $\mu^{-1}(K')<A$. A particular
case is the following. Let $W$ be a Coxeter group and let $A$ be
the corresponding  Artin group. We define \textit{the colored
Artin group}, $CA$, as the kernel of the canonical epimorphism
$\mu:A\longrightarrow W$. The colored Artin groups of spherical
type appear in the theory of arrangements of hyperplanes as
fundamental groups of complements of the so-called Coxeter
arrangements (see \cite{O}).

\bigskip

\begin{figure}[ht]
  \includegraphics{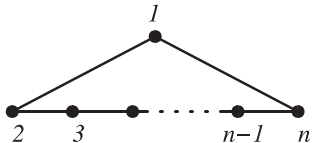}\\
  \caption{Coxeter graph $\tilde{A}_{n-1}$}\label{Fig2}
\end{figure}

\begin{figure}[ht]
  \includegraphics{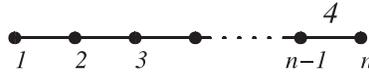}\\
  \caption{Coxeter graph  $B_n$}\label{Fig3}
\end{figure}

\noindent \textbf{Example 3:} The Artin group of type
$\tilde{A}_{n-1}$ is the group associated to the Coxeter graph
pictured in Figure \ref{Fig2},
and the Artin group of type $B_n$ has a  presentation given by the
Coxeter graph pictured in Figure \ref{Fig3}.
The Coxeter group $W(\tilde{A}_{n-1})$ is isomorphic to the
infinite group $\mathbb{Z}^{n-1}\rtimes \Sigma_n$, so the Artin
group $A(\tilde{A}_{n-1})$ is not of spherical type. Let
$C_2=\{\pm 1\}$ denote the cyclic group of order 2. On the other
hand, the Coxeter group of type $B_n$ is isomorphic to the
semi-direct product $(C_2)^n\rtimes \Sigma_n$ which is finite,
hence $A(B_n)$ is of spherical type. Let $\{a_1,\ldots,a_n\}$ be
the standard generators of $A(\tilde{A}_{n-1})$, and let
$\{b_1,\ldots,b_n\}$ be the standard generators of $A(B_n)$.
Consider the homomorphism $\Phi:A(\tilde{A}_{n-1})\longrightarrow
A(B_n)$ defined by
 $$\Phi(a_i)=b_i,\quad i=1,\cdots, n-1$$
 $$\Phi(a_{n})=b_{n}^{-1}b_{n-1}^{-1}\cdots b_{2}^{-1}b_1 \cdots
 b_n.$$
The homomorphism $\Phi$ is injective (see \cite{AL}) allowing us
to see the group $A({\tilde{A}_n})$ as a subgroup of $A(B_n)$.

Define the homomorphism $\varphi:A(B_n)\longrightarrow \mathbb{Z}$
  by $\varphi(b_n)=1$ and $\varphi(b_i)=0$ for $i<n$. In \cite{KP} it is proved that
  the group $A({\tilde{A}_n})$ is the inverse image of $\{1\}$ by $\varphi$. Now, because
$\mathbb{Z}$ is abelian, hence it is endowed with a solution to
the membership problem, and $\{1\}$ is finite, the triple
$(\mathbb{Z},\{1\},\varphi)$ verifies Conditions \textbf{PC}. So,
we can apply our theory to the pair
$(A({\tilde{A}_{n-1}}),A(B_n))$, and give an explicit solution to
the conjugacy problem in $A({\tilde{A}_{n-1}})$. Another solution
to the conjugacy problem in $A({\tilde{A}_{n-1}})$ can be derived
from the facts that $A({\tilde{A}_{n-1}})$ is biautomatic (see
\cite{Ch}), and that a biautomatic structure on a group provides a
solution to the conjugacy problem (see \cite{E}). The centralizers
of some specific elements of $A({\tilde{A}_{n-1}})$ are given in
\cite{Di}.

\vspace{0.15in}

\noindent \textbf{Example 4:} Let $X=\{1,\ldots,m\}\subset
\{1,\ldots,n\}$, and let $\pi_X:\mathcal{B}_n(X)\longrightarrow
C\mathcal{B}_{m}$ be the homomorphism which sends a braid
$b=(b_1,\ldots,b_n)\in \mathcal{B}_n(X)$ to the pure braid
$\pi_X(b)=(b_1,\ldots,b_m)$. It is easily checked that $\pi_X$ is
a well-defined homomorphism. Define $\mathcal{IB}_n(X)$ as the
kernel of the projection $\pi_X$. Note that $\mathcal{IB}_n(X)$ is
the braid group on $n-m$ strands of the $m$ punctured disc.
Presentations for  these groups are given in \cite{L}.

\begin{figure}[ht]
  \includegraphics{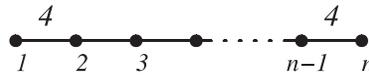}\\
  \caption{Coxeter graph of $\tilde{C}_{n-1}$}\label{Fig4}
\end{figure}

Assume $X=\{1,2\}$. Then the Artin group of type
$\tilde{C}_{n-1}$, defined by the Coxeter graph pictured in Figure
\ref{Fig4},
is isomorphic to $\mathcal{B}_{n}(X)$ (see \cite{AL}). By Example
1, we have that the group $\mathcal{B}_{n}(X)$ has a solution to
the explicit conjugacy problem, and  we can compute a generating
set of the centralizer of $b$ for all $b\in\mathcal{B}_n(X)$.
Moreover, the colored braid group $C\mathcal{B}_2$ is isomorphic
to $\mathbb{Z}$, thus the triple $(C\mathcal{B}_2,\{1\},\pi_X)$
verifies Conditions \textbf{PC}, hence we can apply our theory to
the subgroup $\mathcal{IB}_{n}(X)< \mathcal{B}_{n}(X)$, and obtain
a solution to the explicit conjugacy problem in
$A(\tilde{C}_{n-1})$. As far as I know, there is no other known
solution to the conjugacy problem in $A(\tilde{C}_{n-1})$.

Beside the finite groups, another family of well-understood
Coxeter groups are the so-called affine Coxeter groups which
include, notably, the groups of type $\tilde{A}_n$ and
$\tilde{C}_n$. The corresponding Artin groups, called
\textit{affine type Artin groups}, are poorly understood. Some
other recent results concerning the groups $A(\tilde{A}_n)$ and
$A(\tilde{C}_n)$ can be founded in  \cite{CC}, \cite{Ch} and
\cite{Di}. The previous techniques (Examples 1 and 3) give
algorithms to treat the cases $\tilde{A}_n$ and $\tilde{C}_n$, but
we do not know how to deal with the other affine type Artin
groups.

Now, we show that it is not possible to apply our theory to the
pair $(\mathcal{IB}_n(X),\mathcal{B}_n(X))$ for all $m\geq 5$. We
do not know how to treat the cases $m=3$ and $m=4$.

\smallskip

\noindent\textbf{Remark:} If $m=n$ then $\mathcal{IB}_n(X)$ is the
trivial group, hence we have a solution to the conjugacy problem.
If $m=n-1$ then $\mathcal{IB}_n(X)\cong F_m$ (the free group on
$m$ generators), thus we also have a solution to the conjugacy
problem in $\mathcal{IB}_n(X)$. If $3\leq m\leq n-2$ then we do
not know any solution to the conjugacy problem in
$\mathcal{IB}_n(X)$.

\begin{proposition}\label{proIB}

Let $X=\{1,\ldots,m\}\subset \{1,\ldots,n\}$, $m\geq 5$. There is
no triple $(K,K',\phi)$ such that $K$ is a group endowed with a
solution to the membership problem, $K'$ is a subgroup of $K$,
$\phi$ is a homomorphism $\phi:\mathcal{B}_{n}(X)\longrightarrow
K$,  and $\phi^{-1}(K')=\mathcal{IB}_{n}(X)$.

\end{proposition}

\begin{proof}

Note that, if such a triple exists, then we can assume that the
homomorphism $\phi$ is surjective. If $\phi$ is not surjective,
and $(K,K',\phi)$ verifies Condition \textbf{PC}, then the triple
$(Im(\phi),K',\phi)$ also verifies Condition \textbf{PC} and, in
this case, $\phi$ is surjective.

Consider the natural embedding $\iota:\mathcal{B}_m
\longrightarrow \mathcal{B}_n(X)$ defined by
$\iota(\sigma_i)=\sigma_i$ for all $i=1,\ldots,m-1$. Then $\iota$
restricts to a monomorphism $\iota:C\mathcal{B}_m \longrightarrow
\mathcal{B}_n(X)$ which, actually, is a section of
$\pi_X:\mathcal{B}_n(X) \longrightarrow C\mathcal{B}_m$. So, we
have the split exact sequence,

$$\xymatrix{ {1}\ar@{->}[r] &  {\mathcal{IB}_{n}(X)}\ar@{^{(}->}[r] &
{\mathcal{B}_{n}(X)}\ar@/_/@{->}[r]_{\pi_X}   &
{C\mathcal{B}_{m}}\ar@{->}[r] \ar@/_/@{->}[l]_{\iota}  & 1 }.$$

Define the homomorphism $\pi'_X:K\longrightarrow \mathcal{B}_m$ as
follows. Let $v\in K$. Take $\tilde{v}\in \mathcal{B}_{n}(X)$ such
that $\phi(\tilde{v})=v$, and set $\pi'_X(v)=\pi_X(\tilde{v})$.
Suppose that $\tilde{w}$ is another element of
$\mathcal{B}_{n}(X)$ such that $\phi(\tilde{w})=v.$ Then
$\phi(\tilde{w}\tilde{v}^{-1})=1$, thus, since $1\in K'$,
$\tilde{w}\tilde{v}^{-1}\in \phi^{-1}(K')=\mathcal{IB}_{n}(X)$,
therefore $\pi_X(\tilde{w}\tilde{v}^{-1})=1$, that is,
$\pi_X(\tilde{w})=\pi_X(\tilde{v})$. This shows that the
homomorphism $\pi'_X$ is well-defined. Remark that, by
construction, $\pi'_X\circ \phi=\pi_X.$

Now we prove that $Ker \pi'_X=K'$. Let $a\in K'$. Take
$\tilde{a}\in \mathcal{IB}_{n}(X)$ such that $\phi(\tilde{a})=a$.
So, $\pi'_X(a)=\pi_X(\tilde{a})=1$, hence $a\in Ker \pi'_X$. Let
$a\in Ker \pi'_X$. Then there exists $\tilde{a}\in
\mathcal{B}_{n}(X)$ such that $\phi(\tilde{a})=a$. Because
$\pi_X(\tilde{a})=\pi'_X(a)=1$, we have  $\tilde{a}\in
\mathcal{IB}_{n}(X)$ and $\phi(\tilde{a})=a\in K'$.

Let $\iota'=\phi\circ\iota$, and observe that $\iota'$ is a
section of $\pi'_X$.

At this point, we have proved that the following diagram is
commutative and that the lines are split exact sequences.

$$\xymatrix{ {1}\ar@{->}[r] &  {\mathcal{IB}_{n}(X)}\ar@{^{(}->}[r]\ar@{->}[d]^{\phi} &
{\mathcal{B}_{n}(X)}\ar@/_/@{->}[r]_{\pi_X} \ar@{->>}[d]^{\phi}  &
{C\mathcal{B}_{m}}\ar@{->}[r] \ar@/_/@{->}[l]_{\iota} \ar@{->}[d]^{id} & 1 \\
{1}\ar@{->}[r]  & {K'}\ar@{^{(}->}[r] & {K} \ar@/_/
@{->}[r]_{\pi'_X}
 &
{C\mathcal{B}_{m}}\ar@{->}[r] \ar@/_/@{->}[l]_{\iota'}  & 1 }$$

Let $F_2$ be the free group on two generators,  and let
$R=F_2\times F_2$.  Theorem 1 of \cite{Ma} tells that we can
inject $R$ in $C\mathcal{B}_{m}$, if $m\geq 5$, and Theorem 1 of
\cite{Mi} says that there is no solution to the membership problem
in $R$. Now, $C\mathcal{B}_m$ can be embedded in $K$ via $\iota'$,
thus $R$ can be embedded in $K$, too, so,  $K$ has no  solution to
the membership problem.

\end{proof}

\section{The conjugacy algorithm}

In this section we  assume that $G$ is a group, $H$ is a subgroup
of $G$, and $(K,K',\phi)$ is a triple that verifies Condition
\textbf{PC}, and  we  prove an algorithm which solves the explicit
conjugacy problem in $H$.

We start with a preliminary algorithm.

\bigskip

\textbf{Algorithm $3.1$. From conjugation in $G$ to conjugation in
$H$}

\bigskip
Let $a,b\in H$, $h\in G$ such that $b=h^{-1}ah$. This algorithm
will return FALSE, if $a$ and $b$ are not conjugated in $H$, and
TRUE together with an element $h'\in H$ such that
$b=(h')^{-1}ah'$, otherwise. This algorithm works as follows.

Let $a,b\in H$, and $h\in G$ such that $b=h^{-1}ah$. Recall that
the group $G$ is endowed with an algorithm that computes a finite
generating set $D$ of $C_G(a)$. By Condition \textbf{PC} the set
$\phi(C_G(a)h)\cap K'$ is finite, and we can enumerate all its
elements as follows. Observe that
$\phi(C_G(a))=\langle\phi(D)\rangle$. For each $k\in K'$, check
whether  $kh^{-1}\in \phi(C_G(a))$. This can be done because we
know that $\phi(C_G(a))$ is generated by $\phi(D)$, and because
$K$ is endowed with a solution to the membership problem. Since
$K'$ is finite this process will stop. If $\phi(C_G(a)h)\cap K'=
\varnothing$, then $a$ and $b$ are not conjugated in $H$. Indeed,
if $b=(h')^{-1}ah'$ for some $h'\in H$, then $h'h^{-1}\in C_G(a)$,
and $\phi(h^{\prime})=\phi(h^{\prime}h^{-1}h)\in K'$, so
$\phi(C_G(a)h)\cap K'\neq \varnothing$. Suppose that
$\phi(C_G(a)h)\cap K'\neq \varnothing$. Take an element $p\in
\phi(C_G(a)h)\cap K'$ and write
$p=\phi(d_1)\cdots\phi(d_q)\phi(h)$, where $d_1,\ldots,d_q\in D$.
Set  $d=d_1\cdots d_q$, and $h'=dh$. Then $\phi(h')=p\in K'$ (thus
$h'\in H$), and $(h')^{-1}ah'=h^{-1}d^{-1}adh=h^{-1}ah=b$. We can
resume this algorithm in the following way.

\bigskip

INPUT: A group $G$, a subgroup $H$ of $G$, a triple $(K,K',\phi)$
that verifies Condition \textbf{PC}, two elements $a,b\in H$, an
element $h\in G$  such that $b=h^{-1}ah$.

\bigskip

\begin{enumerate}
    \item Compute a generating set $D$ of $C_G(a).$
    \item For each $p\in K'$ do:
    \begin{enumerate}
        \item If $p\phi(h^{-1})\in \phi(C_G(a))$, then do:
        \begin{enumerate}
             \item Write $p=\phi(d_1)\cdots\phi(d_q)\phi(h)$, where $d_1,\ldots,d_q\in
         D$.
             \item return TRUE and $h'=d_1\cdots d_qh$.
        \end{enumerate}
    \end{enumerate}
    \item return FALSE
\end{enumerate}

\bigskip

OUTPUT: either FALSE if $a$ and $b$ are not conjugated in $H$, or
TRUE together with an element $h'\in H$ such that $b=h'^{-1}ah'$,
otherwise.

\bigskip

\textbf{Algorithm $3.2$. Explicit conjugacy algorithm in $H$}

\bigskip

Given $a,b\in H\subset G$, apply the known conjugacy algorithm to
$a$ and $b$. If $a$ and $b$ are not conjugated as elements in $G$,
then they cannot be conjugated in $H$. Otherwise, take  $h\in G$
such that $h^{-1}ah=b$, and apply Algorithm $3.2$.

\bigskip

INPUT: A group $G$, a subgroup $H$ of $G$, a triple $(K,K',\phi)$
that verifies Condition \textbf{PC}, two elements $a,b\in H$.

\bigskip

\begin{enumerate}
    \item For $a,b\in H$ test if they are conjugated in $G$. If they are not then return FALSE.
    \item Let $h\in G$ be the output of the explicit conjugacy
    algorithm of $G$.
    \item Apply Algorithm $3.1$ to $a,b$ and $h$ and return the output of Algorithm $3.1$
\end{enumerate}

\bigskip

OUTPUT: either FALSE if $a$ and $b$ are not conjugated in $H$, or
TRUE together with an element $h'\in H$ such that $b=h'^{-1}ah'$,
otherwise.

\bigskip

\section{Considerations on the complexity of the algorithms}

In general it is not possible to present an upper bound for the
complexity of the algorithms presented in the previous section.
There are several factors to take in charge. The first ones are
the complexities of the conjugacy problem and of the algorithm to
compute a generating set of the centralizers in $G$. Theses
complexities, in the cases of the Artin groups of spherical type,
are related to the size of a particular subset of the conjugacy
class, the number of generators and the length of the word (see
\cite{FGM}, \cite{FGM2} and \cite{Ge}). Regarding the algorithms
presented, we have also the following factors, the complexity of a
solution to the membership problem in $K$, and the number of
generators of the centralizers in $G$. Regarding this last factor
we have some upper bounds in the case of the braid group (see
\cite{GMW}).


\smallskip\noindent {\bf Nuno Franco},

\smallskip\noindent
R. Rom\~ao Ramalho 59, Departamento de Matem\'atica, CIMA-UE,
Universidade de \'Evora, 7000 \'Evora, Portugal

\smallskip\noindent Institut de Math\'ematiques de Bourgogne, UMR 5584 du
CNRS, Universit\'e de Bourgogne, B.P. 47870, 21078 Dijon cedex,
France

\smallskip\noindent E-mail: {\tt nmf@uevora.pt}

\end{document}